\documentclass[12pt]{amsart}
\usepackage{amscd,amssymb}
\usepackage[arrow,matrix]{xy}
\usepackage[plainpages,backref,urlcolor=blue]{hyperref}

\theoremstyle{plain}
\newtheorem{thm}[subsection]{Theorem}
\newtheorem{lem}[subsection]{Lemma}
\newtheorem{prop}[subsection]{Proposition}
\newtheorem{cor}[subsection]{Corollary}

\theoremstyle{definition}

\newtheorem{definition}[subsection]{Definition}
\newtheorem{ex}[subsection]{Example}
\newtheorem{question}[subsection]{Question}

\numberwithin{equation}{section}
\newcommand{\F}{{\mathcal F}}

\newcommand{\B}{{\mathcal B}}
\newcommand{\QQ}{{\mathcal Q}}

\newcommand{\CC}{{\mathcal C}}

\newcommand{\E}{{\mathcal E}}

\newcommand{\al}{{\alpha}}
\newcommand{\be}{{\beta}}

\newcommand{\R}{\mathbb{R}}
\newcommand{\C}{\mathbb{C}}

\newcommand{\PP}{\mathbb{P}}

\DeclareMathOperator{\im}{im}

\DeclareMathOperator{\mult}{mult}

\DeclareMathOperator{\codim}{codim}

\begin{document}

\title [Chebyshev curves, free resolutions and rational curve arrangements]
{Chebyshev curves, free resolutions and rational curve arrangements }

\author[Alexandru Dimca]{Alexandru Dimca$^1$}
\address{Institut Universitaire de France et Laboratoire J.A. Dieudonn\'e, UMR du CNRS 6621,
                 Universit\'e de Nice Sophia-Antipolis,
                 Parc Valrose,
                 06108 Nice Cedex 02,
                 France}
\email{dimca@unice.fr}

\author[Gabriel Sticlaru]{Gabriel Sticlaru}
\address{Faculty of Mathematics and Informatics,
Ovidius University,
Bd. Mamaia 124, 900527 Constanta,
Romania}
\email{gabrielsticlaru@yahoo.com }
\thanks{$^1$ Partially supported by the French-Romanian Programme LEA Math-Mode
and  ANR-08-BLAN-0317-02 (SEDIGA)}



\begin{abstract} First we construct a free resolution for the Milnor (or Jacobian) algebra $M(f)$ of
a complex projective Chebyshev plane curve $\CC_d:f=0$ of degree $d$. In particular, this resolution implies that the dimensions of the graded components $M(f)_k$ are constant for $k \geq 2d-3.$

Then we show that the Milnor algebra of a nodal plane curve $C$ has such a behaviour if and only if all
the irreducible components of $C$ are rational.

For the Chebyshev curves, all of these components are in addition smooth, hence they are lines or conics and explicit factorizations are given in this case.

\end{abstract}

\maketitle


\section{Introduction and statement of results} \label{sec:intro}

Let $T_d(x)=\cos (d \arccos (x))$ be the Chebyshev polynomial of degree $d$, some of whose properties will be recalled in Section 2 below.

The main of these properties is that $T_d(x)$ has only two critical values, namely $\pm 1$.
This was used by S. V. Chmutov to construct complex projective hypersurfaces with a large number of nodes, i.e. $A_1$-singularities, see \cite{AGV}, volume 2, p. 419 and \cite{Chm}. For another occurence of Chebyshev polynomials in modern algebraic geometry see \cite{BCG}.

In this paper we consider the complex Chebyshev curves $\CC_d$ defined in the complex projective plane $\PP^2$ by the equation $T_d(x,y,z)=0$, where $T_d(x,y,z)$ is the polynomial obtained by homogenization of the polynomial $F^+_d(x,y)=T_d(x)+T_d(y)$, hence exactly the 1-dimensional case of Chmutov's hypersurfaces.

It turns out that these curves are unions of lines and conics in general position, see Corollaries \ref{corA} and \ref{corB}. Hence their study may be of interest in the theory of rational curve arrangements, see for instance \cite{Art}, \cite{Cog}, \cite{ST}. There are interesting relations with knot theory as well, see \cite{KPR}.

Let $S=\C[x,y,z]$ be the graded ring of polynomials in $x,y,z$ with complex coefficients and denote by $S_r$ the vector space of homogeneous polynomials in $S$ of degree $r$.
For any polynomial $f \in S_r$ we define the {\it Jacobian ideal} $J_f \subset S$ as the ideal spanned by the partial derivatives $f_x,f_y,f_z$ of $f$ with respect to $x,y,z$ and the graded Milnor (or Jacobian) algebra
\begin{equation} 
\label{eq1}
M(f)=S/J_f.
\end{equation}
The study of such Milnor algebras is related to the singularities of the corresponding projective curve $C_f:f=0$, see \cite{CD}, as well as to the mixed Hodge theory of 
the curve $C_f$ and of its complement $\PP^2 \setminus C_f$, see the foundational article by Griffiths \cite{Gr} and also \cite{DS}, \cite{DSW}, all three papers treating the case of hypersurfaces in $\PP^n$.

The Hilbert-Poincar\'e series of a graded $S$-module $N$ of finite type is defined by 
\begin{equation} 
\label{eq2}
HP(N)(t)= \sum_{k\geq 0} \dim N_kt^k
\end{equation}
and it is known, see for instance \cite{GP}, to be a rational function of the form 

\begin{equation} 
\label{eq3}
HP(N)(t)=\frac{P(N)(t)}{(1-t)^3}.
\end{equation}
In other words, to determine the series $HP(N)(t)$ it is enough to determine the polynomial $P(N)(t)$.

The best way to understand a Milnor algebra $M(f)$ is to construct a free resolution. Our first main result gives such a resolution for the case $f=T_d(x,y,z)$. The result depends on the parity of $d$.

\begin{thm}
\label{thm1}

\noindent (i) If $d=2m$ is even, then there is a free resolution of the Milnor algebra $M=M(T_d(x,y,z))$
of the form
$$ 0 \to S^m(-4m+1) \to S^{m-1}(-4m+3)\oplus S^3(-4m+2)\to S^3(-2m+1) \to S \to M \to 0.$$
In particular,
$$P(M)(t)=1-3t^{2m-1}+(m-1)t^{4m-3}+3t^{4m-2}-mt^{4m-1}.$$

\medskip

\noindent (ii) If $d=2m+1$ is odd, then there is a free resolution of the Milnor algebra $M=M(T_d(x,y,z))$
of the form
$$ 0 \to S^m(-4m-1) \to S^{m}(-4m+1)\oplus S^2(-4m)\to S^3(-2m) \to S \to M \to 0.$$
In particular,
$$P(M)(t)=1-3t^{2m}+mt^{4m-1}+2t^{4m}-mt^{4m+1}.$$

\end{thm}

Here $S^k=S \oplus...\oplus S$, the direct sum of $k$ copies of $S$, is endowed with the usual grading, namely $(S^k)_r=S_r \oplus...\oplus S_r$ and $(S^k(e))_r=(S^k)_{e+r}$. Moreover, linear morphisms $u:N \to N'$ between graded modules are supposed to preserve the grading, i.e. $u(N_r) \subset N'_r$.

Note that in both cases $\deg P(M)(t)=2d-1$ and $P(M)(t)$ is divisible by $(t-1)^2$. To see this, one 
may use Theorem 5.3.7 in \cite{GP}, since the Krull dimension of $M$ is clearly one, as the corresponding zero-set $V(J_f)$ is a union of lines in $\C^3$.
It follows that in fact one may write 
$$HP(M)=\frac{Q(M)(t)}{1-t}$$
where $Q(M)(t)$ is a polynomial of degree $2d-3$.
Using the relation between Hilbert-Poincar\'e series and Tjurina numbers obtained in \cite{CD}, we get the following.

\begin{cor}
\label{cor1}
\noindent (i) If $d=2m$ is even, then the Chebyshev curve $\CC_d$ has $2m(m-1)$ nodes as singularities and
$\dim M_k=2m(m-1)$
for $k \geq 2d-3$.
\medskip

\noindent (ii) If $d=2m+1$ is odd, then the Chebyshev curve $\CC_d$ has $2m^2$ nodes as singularities and
$\dim M_k=2m^2$
for $k \geq 2d-3$.

\end{cor}

The results proved in \cite{CD} for hypersurfaces with isolated singularities imply that we have stabilization at the step $k_0$, i.e. the dimensions $\dim M_k$ become constant for $k \geq k_0$
with $k_0 =3(d-2)+1$ in general, hence a much weaker result than in the case of Chebyshev curves.

It is surprising that the stabilization at $k_0=2d-3$ for a nodal plane curve has a clear-cut geometric description given in the following, which is the second main result of our note. 

\begin{thm}
\label{thm2}
Let $C \subset \PP^2$ be a nodal curve and $f=0$ be a reduced equation for $C$.

Then the following are equivalent.

\noindent (i) The sequence of dimensions $\dim M(f)_k$ stabilizes at $k_0=2d-3$.

\noindent (ii) $\dim M(f)_{2d-3}=n(C)$, the total number of nodes, alias $A_1$-singularities of the curve $C$. 

\noindent (iii) Any irreducible component $C_i$  of $C$ is a rational curve. In particular, if $C_i$ is smooth, then $C_i$ is either a line or a conic. 

\end{thm}

In other words, $(ii)$ is a numerical test to decide whether a given nodal curve $C$ is an arrangement of rational curves, which is the same as asking that $C$ is covered by a finite number of rational parametrizations.

\bigskip

The proof of Theorem \ref{thm1} is elementary, based on a detailed analysis of the geometry of the Chebyshev curve $\CC_d$ and of an associated curve $\B_d$, defined by the difference of two Chebyshev polynomials. On the other hand, the proof of Theorem \ref{thm2} uses in two points subtle facts of mixed Hodge theory, once in the proof of Lemma \ref{lem1} and then at the end of the proof, see the formula \eqref{eqDSW}. The study of the relations between the Hilbert-Poincar\'e series of the Milnor algebra $M(f)$ and Hodge theory is continued in our recent preprint \cite{DSt}.

Numerical experiments with the CoCoA package \cite{Co} and the Singular package \cite{Sing} have played a key role in the completion of this work.

\section{Basic facts on Chebyshev polynomials} \label{sec:two}

The $d$-th Chebyshev polynomial $T_d(x)$ has $d-1$ critical points, namely $\lambda_k=\cos (k\pi/d)$ for $k=1,...,d-1$. Hence $T_d'(\lambda_k)=0$ and $T_d(\lambda_k)=(-1)^k$. Note also that $T_d'' (\lambda_k) \ne 0$ for $d \geq 2$.

It follows that, for $d=2m+1$ odd,  the critical values $\pm 1$ are both attained $m$ times.
When $d=2m$, the critical value $ 1$ is attained $m-1$ times, while $-1$ is attained $m$ times.

The Chebyshev curve $\CC_d$ has as singular points exactly the points $(\lambda_p:\lambda_q:1) \in \PP^2$
such that $p+q$ is odd. The above remarks show that all these singularities are nodes (i.e. singularities of type $A_1$) and their total number is as stated in Corollary \ref{cor1}.

Note that all the singularities of the Chebyshev curve $\CC_d$ are in fact situated in the affine plane $\C^2 \subset \PP^2$ given by $z=1$. We use $x,y$ as coordinates on this affine plane and freely identify a homogeneous polynomial $f \in S_r$ to a polynomial $g(x,y)=f(x,y,1) \in \C[x,y]_r$,
the set of polynomials in $\C[x,y]$ of degree at most $r$.

The singularities of $\CC_d$ corresponds to the set of points in $\C^2$ given by
\begin{equation} 
\label{eq4}
A_d=\{(\lambda_p,\lambda_q)~~:~~0<p<d,~~0<q<d, ~~p+q \text{ odd } \}.
\end{equation}
We consider also the complementary set
\begin{equation} 
\label{eq5}
B_d=\{(\lambda_p,\lambda_q)~~:~~0<p<d,~~0<q<d, ~~p+q \text{ even } \}.
\end{equation}
Note that the affine curve $\B_d: F_d^-(x,y):=T_d(x)-T_d(y)$ has only nodes as singularities and they are located exactly at the points in $B_d$. The polynomial $F_d^-(x,y)$ has a wonderful factorization,
established in \cite{KPR}. 

\begin{prop}
\label{prop1}
The polynomial $F^-_d(x,y)$ is a product of linear and quadratic factors. More precisely, if we set
$$g_k(x,y)= x^2-2\lambda_{2k}xy+y^2-(1-\lambda_{2k}^2),$$
then one has the following.

\noindent (i) If $d=2m$ is even, then the linear factors are $x-y$ and $x+y$, and the irreducible quadratic factors are $g_k(x,y)$ for $k=1,...,m-1$.

\medskip

\noindent (ii) If $d=2m+1$ is odd, then  the only linear factor is $x-y$, and the irreducible quadratic factors are $g_k(x,y)$ for $k=1,...,m$.

\end{prop}

Recall that $T_d(x)$ is an even (resp. odd) polynomial when $d$ is even (resp. odd). Hence, for $d$ odd, we get
$$F_d^+(x,y)=F^-_d(x,-y).$$
Using Proposition \ref{prop1} we get the following.

\begin{cor}
\label{corA}
If $d=2m+1$ is odd, then the Chebyshev polynomial $F_d^+(x,y)$ splits as a product of the linear factor $x+y$ and the quadratic factors $g^+_k(x,y)= x^2+2\lambda_{2k}xy+y^2-(1-\lambda_{2k}^2),$ for $k=1,...,m.$
In particular, the polynomials $F_d^+(x,y)$ and $F^-_d(x,y)$ are affinely equivalent.
\end{cor}

The same method of proof as for Proposition \ref{prop1} yields the following.

\begin{cor}
\label{corB}
If $d=2m$ is even, then the Chebyshev polynomial $F_d^+(x,y)$ splits as a product of the quadratic factors $g^+_k(x,y)= x^2+2\lambda_{2k-1}xy+y^2-(1-\lambda_{2k-1}^2),$ for $k=1,...,m.$
In particular, in this case the polynomials $F_d^+(x,y)$ and $F^-_d(x,y)$ are not affinely equivalent.
\end{cor}

Note that each quadratic form $g_k(x,y)$ defines an ellipse $\E_k$ in the real plane $\R^2$, whose symmetry axes are the lines $L_1:x-y=0$ and $L_2:x+y=0$. In particular, the intersections $L_j \cap \E_k$ consists of two distinct points in $B_d$.

Moreover, two distinct ellipses $\E_k$ and $\E_{\ell}$ intersect in precisely four distinct points given by 
\begin{equation} 
\label{eq6}
(\lambda_{k+\ell}, \lambda_{k-\ell}), ~~ (-\lambda_{k+\ell}, -\lambda_{k-\ell})~~, (\lambda_{k-\ell}, \lambda_{k+\ell}),\text{ and }
(-\lambda_{k-\ell}, -\lambda_{k+\ell}).
\end{equation}
Notice that this $4$ points determine a rectangle, call it $R_{k,\ell}$.

\section{Construction of the free resolutions} \label{sec:3}

Let $\F(B_d)$ be the vector space of $\C$-valued functions defined on $B_d$.

For each positive integer $r$, we have an evaluation map
\begin{equation} 
\label{eq7}
e(d,r): \C[x,y]_r \to \F(B_d)
\end{equation}
obtained by sending a polynomial $h \in \C[x,y]_r$ to the function $x \mapsto h(x)$ for $x \in B_d$.
We have the following interpolation-type result, similar in spirit to \cite{Xu}.

\begin{prop}
\label{prop2}
The map $e(d,r)$ is injective for $r \leq d-3$ and surjective for $r \geq d-2$.

\end{prop}

\proof

To prove the injectivity part, it is enough to treat the case $r=d-3$. Let $Q \in \C[x,y]_{d-3}$ be a nonzero polynomial vanishing at any point in $B_d$ and let $\QQ$ be the corresponding complex plane curve. If $L$ is a line which is an irreducible component of the curve $\B_d$, then $L$ contains exactly $d-1$ points in $B_d$. It follows that
$$\sum_{x \in L\cap \QQ}\mult_x(L, \QQ) \geq d-1>\deg L \cdot \deg Q=d-3.$$
It follows by Bezout's Theorem, see for instance \cite{GH}, p. 172, that $L$ is an irreducible component of $\QQ$.

Consider next the conic (which is the complex version of the ellipse) $\E_k$. It meets any other conic $\E_{\ell}$ for $\ell \ne k$ in four points. 

In the case $d=2m$ even, we get in this way $4(m-2)$ points of this type on $\E_k$, all of them in $B_d$. There are four more points coming from the intersections of $\E_k$ with the lines $L_1$ and $L_2$. In all we get $4(m-1)$
points of intersection between $\E_k$ and $\QQ$. One has
$$\sum_{x \in \E_k \cap \QQ}\mult_x(\E_k, \QQ) \geq 4(m-1)>\deg \E_k \cdot \deg Q=2(2m-3).$$
It follows by Bezout's Theorem that $\E_k$ is an irreducible component of $\QQ$.

In the case $d=2m+1$ even, we get in this way $4(m-1)$ points of intersection with the other conics on $\E_k$. There are $2$ more points coming from the intersections of $\E_k$ with the line $L_1$. In all we get $4m-2$
points of intersection between $\E_k$ and $\QQ$. One has
$$\sum_{x \in \E_k \cap \QQ}\mult_x(\E_k, \QQ) \geq 4m-2>\deg \E_k \cdot \deg Q=2(2m-2).$$
It follows again by Bezout's Theorem that $\E_k$ is an irreducible component of $\QQ$.

The above shows that $Q$ is divisible by the polynomial $F^-_d(x,y)$, which is a contradiction, since
$\deg Q < \deg F^-_d(x,y)$.

\medskip

To prove the surjectivity, it is enough to suppose $r=d-2$. For each point $b \in B_d$ we construct a polynomial $h_b \in \C[x,y]_{d-2}$ such that $h_b(b) \ne 0$ and $h_b(c)=0$ for any $c \in B_d$, $c \ne b$.

Since all the singularities of $\B_d$ are nodes, it follows that at any point $b \in B_d$ there are exactly two irreducible components of $\B_d$ meeting there. It follows that there are three cases to consider.

\medskip

\noindent (i) the point $b$ is the intersection of two lines (then of course $b=(0,0)$, $d=2m$ and the lines are $L_1$ and $L_2$). Define in this case
$$h_b(x,y)=\frac{F^-_d(x,y)}{x^2-y^2}.$$
\medskip

\noindent (ii) the point $b$ is the intersection of one line $L$ (with equation $\ell(x,y)=0$)
and a conic $\E_k$. Let $b'$ be the second intersection point of these two curves and let $\ell'(x,y)=0$ be the equation of a line passing through $b'$ but different from $L$. Then one can take
$$h_b(x,y)=\frac{F^-_d(x,y)\ell'(x,y)}{\ell(x,y) g_k(x,y)}.$$

\medskip

\noindent (iii) the point $b$ is the intersection of two conics $\E_k$ and $\E_{\ell}$. Let $\ell_1(x,y)=0$ and $\ell_2(x,y)=0$ be the equation of the two sides in the rectangle $R_{k, \ell}$
not passing through $b$. Then one can take
$$h_b(x,y)=\frac{F^-_d(x,y)\ell_1(x,y)\ell_2(x,y)}{g_k(x,y)  g_{\ell}(x,y)}.$$

\endproof

Now we prove Theorem \ref{thm1} by constructing rather explicitly the free resolution of the Milnor algebra $M$. We set $f=T_d(x,y,z)$ and use the notations $F_d^-(x,y,z)$, $g_j(x,y,z)$ a.s.o. to denote the homogenized polynomials associated to  $F_d^-(x,y)$, $g_j(x,y)$ a.s.o.

We consider first the case $d=2m$ even. Then we have to construct a resolution
\begin{equation} 
\label{res}
0 \to R_3 \to R_2 \to R_1 \to R_0 \to M \to 0
\end{equation}
where $R_0=S$, $R_1=S^3(-2m+1)$, $R_2=S^{m-1}(-4m+3)\oplus S^3(-4m+2)$ and $R_3= S^m(-4m+1)$.
The first step is the same as in the Koszul complex, namely we define $u_0:R_1 \to R_0$ to be the map
$(a_1,a_2,a_3) \mapsto a_1f_x+a_2f_y+a_3 f_z$.

Then $K_1=\ker u_0$ is the module of linear relations (syzygies) involving the partial derivatives 
$f_x,f_y,f_z$. Among these, there are three trivial relations, coming also from the Koszul complex, namely
$$r_1=(f_y,-f_x,0),~~r_2=(f_z,0,-f_x),~~r_3=(0,f_z,-f_y).$$
For $(c_1,c_2,c_3) \in S^3(-4m+2) $, we define $u_1(c_1,c_2,c_3)=c_1r_1+c_2r_2+c_3r_3$.

Now we look for nontrivial relations, i.e. relations not lying in the submodule $u_1 (S^3(-4m+2))$.
Note that any relation
\begin{equation} 
\label{eq7.5}
a_1f_x+a_2f_y+a_3 f_z=0
\end{equation}
with all $a_j$ homogeneous polynomials in some $S_r$, implies that $a_3$ vanishes on the set $B_d$.
Indeed, $f_x$ and $f_y$ vanishes on the set $A_d \cup B_d$, and $f_z$ is not vanishing at any point in $B_d$.

Let $E(d,r)=\ker e(d,r)$ and recall that it may be thought of as a subspace in $S_r$ by homogenization. Then Proposition \ref{prop2} implies that $E(d,r)=0$ and hence $a_3=0$ if $r <d-2$. But this implies that $a_1=a_2=0$, hence there are no relations in this range.

Consider now the case $r=d-2$. Using Proposition \ref{prop2}, it follows that the vector space $E(d,r)$  has  dimension given by 
$$ \dim \C[x,y]_{d-2} - |B_d|=m(2m-1)-(2m(m-1)+1)=m-1.$$
On the other hand we have obviously the following linear independent elements in $E(d,d-2)$:
\begin{equation} 
\label{eq8}
\al_3^j=\frac{F_d^-(x,y,z)}{g_j(x,y,z)} \text{ for } j=1,...,m-1.
\end{equation}
Since $\al_3^jf_z$ vanishes on $A_d\cup B_d$, using Noether AF+BG Theorem, see for instance \cite{GH}, p. 703, we get the existence
of unique polynomials $\al_1^j$ and $\al_2^j$ in $S_{d-2}$ such that
$$\al_1^jf_x+\al_2^jf_y+\al_3^j f_z=0$$
for $j=1,...,m-1.$
Hence we've got $m-1$ new relations $\rho_j=(\al_1^j,\al_2^j,\al_3^j) \in K_1$ and we use them to define $u_1(b_1,...,b_{m-1})=b_1\rho_1+....+ b_{m-1}\rho_{m-1}$ for any $(b_1,...,b_{m-1}) \in S^{m-1}(-4m+3)$.

In this way the morphism $u_1$ is completely defined and we have to show that any relation in $K_1$ is in fact in the image of $u_1$. So consider a relation as in \eqref{eq7.5}.

\medskip

The first case to check is when $r=d-1$. As above, $a_3 \in E(d,d-1)$, and this vector space has dimension $(m-1)+2m$, 
again by Proposition \ref{prop2}.

Note that $E'=xE(d,d-2)\oplus yE(d,d-2) \oplus zE(d,d-2) \subset E(d,d-1)$ is a $(3m-3)$-dimensional subspace, hence it is enough to find 2 new linearly independent elements.
Recall that $f_x=F_d^-(x,y,z)_x$ and $f_y=-F_d^-(x,y,z)_y$. 

On the other hand we have $F_d^-(x,y,z)=(x-y)(x+y)\prod_{j=1,m-1}g_j(x,y,z)$. It follows that 
$$f_x=F_d^-(x,y,z)_x=F_d^-(x,y,z)( \frac{1}{x-y}+ \frac{1}{x+y} +\sum_{j=1,m-1}\frac{g_j(x,y,z)_x}{g_j(x,y,z)} )$$
and 
$$-f_y=F_d^-(x,y,z)_y=F_d^-(x,y,z)( -\frac{1}{x-y}+ \frac{1}{x+y} +\sum_{j=1,m-1}\frac{g_j(x,y,z)_y}{g_j(x,y,z)} ).$$
It follows that no nontrivial linear combination of $f_x$ and $f_y$ belongs to $E'$, in other words we may use $xg_j$, $yg_j$, $zg_j$, $f_x$ and $f_y$ to get a basis of $E(d,d-1)$.

This shows that starting with any relation in $(K_1)_{d-1}$, we can modify it by a relation in the image of $u_1$ (in fact without using the relation $r_1$) in order to get a relation with $a_3=0$. But then we get $f_x|a_2$ and $f_y|a_1$, i.e. we get a multiple of the relation $r_1$.

\medskip

The next case to investigate is $r=d$.  Again, $a_3 \in E(d,d)$, and this vector space has dimension $(m-1)+2m+2m+1$, 
again by Proposition \ref{prop2}.

We  consider as above the vector space $E''$ spanned by all the product $q\al_3^j$, where $q$ runs through a basis of $S_2$.
But now these products are not linearly independent, in fact by \eqref{eq8}
one clearly has 
$$g_i(x,y,z)\al_3^i=g_j(x,y,z)\al_3^j$$
for any $i,j=1,...,m-1$.
It follows that the relation
$$g_1(x,y,z)\rho_1-g_j(x,y,z)\rho_j$$ 
for $j=2,...,m-1$ has a trivial last component, and as such is a multiple of the relation $r_1$.
It follows that there are homogeneous polynomials $\be_j \in S_1$ such that one has
$$g_1(x,y,z)\rho_1-g_j(x,y,z)\rho_j-\be_jr_1=0,$$
for $j=2,...,m-1$.

Hence we get $m-2$ relations among the fundamental relations $\rho_j$ and $r_k$.
We use them to define the morphism $u_2:R_3\to R_2$ on the first $(m-2)$-components of 
$R_3= S^m(-4m+1)$ by setting 
$$u_2(e_j)=(g_1(x,y,z),..., -g_j(x,y,z),....,-\be_j,0,0)$$
where $e_j$ is the canonical basis of $ S^m(-4m+1)$, $-g_j(x,y,z)$ is placed on the $j$-th component in $S^{m-1}(-4m+3)$ for $j=2,...,m-1$, and $-\be_j$ is placed on the first component of $S^3(-4m+2)$,
and moreover $...$ stand for zero components.
There are two more new relations among the fundamental relations $\rho_j$ and $r_k$.
Using the formulas for $f_x$ and $f_y$ given above, we see that
$$(x+y)(f_x-f_y)=2g_1\al_3^1+ \sum_{j=1,m-1}(x+y)((g_j)_x+(g_j)_y)\al_3^j$$
and
$$(x-y)(f_x+f_y)=2g_1\al_3^1+ \sum_{j=1,m-1}(x-y)((g_j)_x-(g_j)_y)\al_3^j),$$
where we set $g_j=g_j(x,y,z)$, $(g_j)_x=g_j(x,y,z)$ a.s.o. to simplify the notation. 
As above, these relations produce relations among relations:
$$(x+y)(r_2-r_3)+2g_1\rho_1+ \sum_{j=1,m-1} (x+y)((g_j)_x+(g_j)_y)  \rho_j-\gamma_1r_1=0$$
and
$$(x-y)(r_2+r_3)+2g_1\rho_1+ \sum_{j=1,m-1} (x-y)((g_j)_x-(g_j)_y)  \rho_j-\gamma_2r_1=0$$
with $\gamma_k \in S_1$.
Using them, we complete the definition of $u_2$ by setting
$u_2(e_{m-1})$ to be 
$$(2g_1+(x+y)((g_1)_x+(g_1)_y),...., (x+y)((g_{m-1})_x+(g_{m-1})_y)  , -\gamma_1,x+y,-(x+y))$$
and
$$u_2(e_{m})=( 2g_1+(x-y)((g_1)_x-(g_1)_y),...., (x-y)((g_{m-1})_x-(g_{m-1})_y)  ,    -\gamma_2,x-y,x-y).$$
The injectivity of $u_2$ follows from the fact that the corresponding matrix
(regarded as a matrix over the field of rational functions $\C(x,y,z)$) has rank $m$.
Indeed, a nonzero size $m$ minor can be obtained by deleting the first and the $m$-th rows in this
matrix.

Returning to the inclusion $(K_1)_d \subset \im (u_1)_d$, the above discussion shows that $\dim E''=3(m-1)-(m-2)=2m-1$
and that it is enough to add the vector space spanned by $xf_x$ and $yf_x$ to $E''$ to get $E(d,d) $.
We conclude that $(K_1)_d \subset \im (u_1)_d$ as above, by first showing that we may assume $a_3=0$.

We leave to the reader to check in the same way the inclusion $(K_1)_r \subset \im (u_1)_r$, for $r>d$ as well as the inclusion $\ker (u_1) \subset \im(u_2)$.

\bigskip

The proof in the case $d=2m+1$ follows the same pattern, with just one major difference.
To construct the required resolution in this case, namely
\begin{equation} 
\label{r1}
R_*: 0 \to S^m(-4m-1) \to S^{m}(-4m+1)\oplus S^2(-4m)\to S^3(-2m) \to S 
\end{equation}
is the same as constructing first a larger resolution
$$R'_*: 0 \to S^m(-4m-1)\oplus S(-4m) \to S^{m}(-4m+1)\oplus S^3(-4m)\to S^3(-2m) \to S $$
and them simplifying it. 

We construct the larger resolution $R'_*$, because its construction is more natural. To start with, the morphism $u_0:S^3(-2m) \to S $ is the same as above, coming from the Koszul resolution.

Let $r=d-2=2m-1$. Then $E(d,r)$ is an $m$-dimensional vector space with a basis given by $\al_3^j$
defined exactly as above, with $j=1,...,m$. They yield $m$ new relations $\rho_j$ with $j=1,...,m$
which are used, together with the trivial relations $r_k$ introduced above, to define the morphism $u_1:S^{m}(-4m+1)\oplus S^3(-4m)\to S^3(-2m)$.

\medskip

Let $r=d-1=2m$, $E(d,r)$ is a $(3m+1)$-dimensional vector space and the vector subspace $E'$ constructed as above has dimension $3m$.

Note that in this case $F_d^-(x,y,z)=(x-y)\prod_{j=1,m}g_j(x,y,z)$. It follows that 
$$f_x=F_d^-(x,y,z)_x=F_d^-(x,y,z)( \frac{1}{x-y} +\sum_{j=1,m}\frac{g_j(x,y,z)_x}{g_j(x,y,z)} )$$
and 
$$-f_y=F_d^-(x,y,z)_y=F_d^-(x,y,z)( -\frac{1}{x-y} +\sum_{j=1,m}\frac{g_j(x,y,z)_y}{g_j(x,y,z)} ).$$
Hence the {\it key new fact in the case $d$ odd} is the relation
$$f_x-f_y=\sum_{j=1,m}(g_j(x,y,z)_x+g_j(x,y,z)_y)\al_3^j.$$
This gives a relation
$$r_2-r_3+\sum_{j=1,m}(g_j(x,y,z)_x+g_j(x,y,z)_y)\rho_j-\delta_1r_1=0$$
which is used to define a morphism $ S(-4m) \to S^{m}(-4m+1)\oplus S^3(-4m)$ which will describe the action of $u_2$ on the factor $S(-4m)$.

For the remaining part of the morphism $u_2$, namely a morphism $S^m(-4m-1) \to S^{m}(-4m+1)\oplus S^3(-4m)$, we need $m$ relations among the relations $\rho_j, r_k$.
The first $(m-1)$ relations are obtained exactly as above, and they have the form
$$g_1(x,y,z)\rho_1-g_j(x,y,z)\rho_j-\be_jr_1=0,$$
for $j=2,...,m$.
The last relation comes from the equality
$$(x-y)f_x=g_1\al_3^1+\sum_{j=1,m}(g_j)_x\al_3^j$$
and hence has the form
$$(x-y)r_2+g_1\rho_1+\sum_{j=1,m}(g_j)_x\rho_j-\delta_2r_1=0.$$

\medskip

To get the formulas for the Hilbert-Poincar\'e series, we start with the resolution \eqref{res} or \eqref{r1} and get
$$HP(M)(t)=HP(R_0)(t)-HP(R_1)(t)+HP(R_2)(t)-HP(R_3)(t).$$
Then we use the well-known formulas $HP(N\oplus N')(t)=HP(N)(t)+HP( N')(t)$, $HP(N(-r))(t)=t^rHP(N)(t)$
and
$$HP(S)(t)=\frac{1}{(1-t)^3}.$$

\section{On rational curve arrangements} \label{sec:4}

The proof of Theorem \ref{thm2} depends on the following standard general fact.
For general properties of mixed Hodge theory we refer to \cite{PS}, see also the corresponding Appendix in \cite{D1} for a brief introduction.

\begin{prop}
\label{prop3}
Let $C \subset \PP^2$ be a nodal curve and set $U=\PP^2\setminus C$. Let $C=\cup_{j=1,r}C_j$ be the decomposition of $C$ as a union of irreducible components, let $\nu_j:\tilde C_j \to C_j$ be the normalization mappings and set $g_j=g(\tilde C_j)$.
Then one has 
$$\dim Gr _F^1H^2(U,\C)=\sum_{j=1,r}g_j$$
and
$$\dim Gr _F^2H^2(U,\C)=\frac{(d-1)(d-2)}{2}.$$

\end{prop}

The second formula follows also from  \cite{DSW}, Theorem (2.2), but we give here a more elementary
independent proof.

\proof

Assume that the curve $C_j:f_j=0$ has degree $d_j$ and has $n_j$ nodes. Recall the definition of the Hodge-Deligne polynomial of a quasi-projective complex variety $X$
$$P(X)(u,v)=\sum_{p,q}E^{p,q}(X)u^pv^q$$
where $E^{p,q}(X)=\sum_s(-1)^s \dim Gr_F^pGr^W_{p+q}H^s_c(X,\C)$, and the fact that it is additive with respect to constructible partitions, i.e. $P(X)=P(X \setminus Y)+P(Y)$, for a Zariski closed subvariety $Y$ in $X$.

Using the normalization maps $\nu_j$, it follows that
$$P(C_j)=P(C_j \setminus (C_j)_{sing})+ P((C_j)_{sing})=P(\tilde C_j \setminus \{2n_j \text { points }\})+n_j=$$
$$=P(\tilde C_j)-P( \{2n_j \text { points }\})+n_j=uv-g_ju-g_jv+1-n_j.$$
Using the fact that the only intersection points are nodes, we get
$$P(C)=P(C_1 \cup....\cup C_r)=\sum_{j=1,r}P(C_j) - \sum_{1 \leq i<j\leq r} P(C_i \cap C_j)=$$
$$
=ruv-(\sum_{j=1,r}g_j)u-(\sum_{j=1,r}g_j)v+r -(\sum_{j=1,r}n_j)- \sum_{1 \leq i<j\leq r}d_id_j.$$
Next we have again by the additivity
$P(U)=P(\PP^2)-P(C)$ where $P(\PP^2)=u^2v^2+uv +1$.
Now, let's look at the cohomology of the smooth surface $U$. The group $H^4_c(U,\C)$, is dual to the group $H^0(U,\C)$, which is $1$-dimensional of type $(0,0)$. It follows that the contribution of 
$H^4_c(U,\C)$ to $P(U)$ is exactly the term $u^2v^2$.
 
The group $H^3_c(U,\C)$, is dual to the group $H^1(U,\C)$, which is $(r-1)$-dimensional of type $(1,1)$. It follows that the contribution of 
$H^3_c(U,\C)$ to $P(U)$ is exactly the term $-(r-1)uv$.
The remaining terms come from the group $H^2_c(U,\C)$, which is dual to the group $H^2(U,\C)$, and this implies that $H^2(U,\C)$
has only classes of type $(2,1)$, $(1,2)$ and $(2,2)$. The dimension $\dim Gr _F^1H^2(U,\C)$
is the number of independent classes of type $(1,2)$ (which correspond to terms in $u$ in $P(U)$), and this gives the first equality.

To establish the second equality we have to work a little more. The dimension $\dim Gr _F^2H^2(U,\C)$
is the number of independent classes of type $(2,1)$ or $(2,2)$ (which correspond to terms in $v$ or to the constant terms in the polynomial $P(U)$). This yields
\begin{equation} 
\label{future}
\dim Gr _F^2H^2(U,\C)=\sum_{j=1,r}(g_j+n_j-1)+\sum_{1 \leq i<j\leq r}d_id_j+1.
\end{equation}
To complete the proof, recall the formula
$$g_j+n_j=p_a(C_j)=\frac{(d_j-1)(d_j-2)}{2},$$
where $p_a$ denotes the arithmetic genus, see \cite{Ha}, p. 298 and p. 54. The result follows
using the relation $d=\sum_{j=1,r}d_j$.

\endproof

We need another preliminary result. Consider the evaluation map
$$e(C)_k:S_k \to \F(\Sigma(C))$$
where $\Sigma(C) \subset \C^3 \setminus \{0\}$ is a finite set in bijection to the nodes of $C$ under the canonical projection $\C^3 \setminus \{0\} \to \PP^2$, $\F(\Sigma(C))$ is the vector space of complex valued functions defined on $\Sigma(C)$, and $e(C)_k(h)$ is the function $s \mapsto h(s)$ for
$h \in S_k$ and $s \in \Sigma(C)$.
\begin{lem} 
\label{lem1}
The mapping $e(C)_{2d-3}$ is surjective.
\end{lem}

\proof
As in the proof of Proposition \ref{prop2} above, it is enough to construct for each point $a \in \Sigma(C)$ a homogeneous polynomial $h_a \in S_p$ for some $p \leq 2d-3$ such that $h_a(a) \ne 0$
and $h_a(b)=0$ for any $b \in \Sigma(C)$, $b \ne a$.
There are two cases to discuss.

\medskip

(i) $a \in C_i \cap C_j$ for $i \ne j$. Then the intersection $C_i\cap C_j$ consists of exactly $d_id_j$ points and a function $h \in S_q$ with $q=d_i+d_j-1$
is in the ideal $I(C_i \cap C_j)$, i.e. $h$ vanishes on this set $C_i\cap C_j$ if and only if can be written (in an unique way) as a sum $h=h_if_i+h_jf_j$ with $\deg(h_i)=d_j-1$ and 
$\deg(h_j)=d_i-1$. It follows that 
$$\dim (S/I(C_i\cap C_j))_q=d_id_j$$
which shows that the corresponding evaluation map $e_q:S_q \to \F(C_i\cap C_j)$ is surjective.
Hence we can find $h' \in S_q$ such that $h'(a) \ne 0$ and $h'(b)=0$ for all $b \in (C_i\cap C_j)$,
$b \ne a$.
Since $a$ is a smooth point on $C_i$ (resp. $C_j$), there is a partial derivative of $f_i$ (resp. of $f_j$) not vanishing at $a$. Denote these partial derivatives by $g_1$ and $g_2$ respectively.

Then the product
$$h_a=h'g_1g_2 \prod_{k=1,r;~~k\ne i;~~k\ne j}f_k$$
has the required properties and $\deg h_a =d+d_i+d_j -3 \leq 2d-3.$

\medskip

(ii) $a$ is a node on the irreducible curve $C_i$. Using the exact sequences (3.13) on p. 201,
(3.16) on p. 202 (for $n=s=2$) and the Example (3.18) on p. 203 (for $m=1$) in \cite{D1}, it follows that the evaluation map
$$e(C_i):S_{d_i-3} \to \Sigma(C_i)$$
can be identified to the surjective map
$$F^2H^2(U,\C)= P^2H^2(U,\C) \to \oplus_kF^2H^2(U_k \setminus C_i,\C)$$
where $P$ denotes the polar filtration on cohomology and $U_k$ a small disc around each singular point of the curve
$C_i$.

Hence it exists a polynomial $h' \in S_{d_i-3}$ such that  $h'(a) \ne 0$ and $h'(b)=0$ for all other nodes $b$ of the curve $C_i$.
Then the product
$$h_a=h'\prod_{k=1,r;~~k\ne i}f_k$$
has the required properties and $\deg h_a =d-3 \leq 2d-3.$ 
\endproof

Next we recall some notation and a key result from \cite{DSW} in our special case.

Let $I \subset S$ be the graded ideal of polynomials vanishing at all the singular points of $C$.
In fact $I_k$ is exactly the kernel of the evaluation map $e(C)_k$, for all $k \geq 0$.
By definition, note that $J_f \subset I$ and also $\dim (S_k/I_k) \leq n(C)$ for all $k$.
Moreover, we have $\dim (S_k/I_k)=n(C)$ exactly when the corresponding evaluation map
$$e(C)_k:S_k \to \F(\Sigma (C))$$
is surjective. The above Lemma shows that this holds for $k=2d-3$.

It is shown in \cite{DSW}, Theorem (2.2) and subsection (2.3), that for a nodal curve in $\PP^2$
one has a natural isomorphism
\begin{equation} 
\label{eqDSW}
Gr _F^1H^2(U,\C)=(I/J_f)_{2d-3}.
\end{equation}
Moreover, it follows from \cite{CD}, Corollary 8, that there are epimorphisms
\begin{equation} 
\label{eqCD}
M(f)_{q-1} \to M(f)_q
\end{equation}
for $q \geq 2d-3$. In particular, $\codim (J_f)_{2d-3}\geq n(C)$, where $\codim$ refers to codimension with respect to $S_{2d-3}.$
This last fact implies that the claim $(i)$ is equivalent to the claim $(ii)$,
since we know that $M(f)_k=\tau(C)=n(C)$ for $k$ large.

If we assume $(ii)$, then by Lemma \ref{lem1} we have  
$$\codim I_{2d-3} = n(C)= \codim (J_f)_{2d-3},$$
and hence we get $(iii)$.

On the other hand, the condition $(iii)$ is saying that $\codim (J_f)_{2d-3}=\codim I_{2d-3}$
which in view of Lemma \ref{lem1} implies $(ii)$.

\end{document}